\documentclass[11pt,reqno]{amsart}
\usepackage[dvips]{epsfig}
\RequirePackage{amsmath}
\usepackage{epsfig}
\usepackage{graphicx,epsfig}
\usepackage{amscd,amsthm,amsfonts,amsopn,amssymb,verbatim}
\usepackage{color}
\usepackage{times}
\numberwithin{equation}{section}
\newtheorem{theorem}{Theorem} 

\newtheorem{lemma}[theorem]{Lemma}

\theoremstyle{remark}

\DeclareMathOperator{\supp}{supp\,}

\def\be{\begin{equation}}
\def\ee{\end{equation}}

\def\ve{\varepsilon}

\allowdisplaybreaks
\def\be{\begin{equation}}
\def\ee{\end{equation}}

\def\ve{\varepsilon}
\def\ve{\varepsilon}

\def\mod{\text{mod\,}}
\allowdisplaybreaks
\parskip =\medskipamount

\begin{document}

\title
{On quadratic irrationals with bounded partial quotients}
\author
{J.~Bourgain}
\address
{School of Mathematics, Institute for Advanced Study, Princeton, NJ 08540}
\address
{bourgain@math.ias.edu}
\thanks{Research supported in part by NSF Grants DMS 1301619}
\begin{abstract}
It is shown that for some explicit constants $c>0, A>0$, the asymptotic for the number of positive 
non-square discriminants $D<x$ with fundamental solution $\ve_D<x^{\frac 12+\alpha}$, $0<\alpha<c$,
remains preserved if we require moreover
$\mathbb Q(\sqrt D)$ to contain an irrational with partial quotients bounded by $A$.
\end{abstract}
\maketitle

\section {Introduction}

We prove the following, related to McMullen's conjecture on irrationals with bounded partial quotients in quadratic number fields.

\noindent
{\bf Theorem 1.}
{\sl For $\alpha_0<\alpha<c$, there are at least $\big(1-o(1)\big)\frac {4\alpha^2}{\pi^2}\sqrt x(\log x)^2$
positive non-square integers $D<x$ for which $\mathbb Q(\sqrt D)$ contains an irrational which partial quotients are bounded by $A$ and such
that moreover the fundamental solution $\ve_D$ to the Pell equation $t^2 -Du^2=1$ is bounded by $x^{\frac 12+\alpha}$.
}

Here $A, c>0$ are explicit constant, and $\alpha_0>0$ arbitrarily small and fixed.
We will explain below the emphasis on the explicitness of these constants.
Theorem 1 certainly holds for $\alpha <\frac 1{200}$.

\noindent
{\bf Remark 1.}
According to the work of Hooley \cite{H} and Fouvry \cite {F}, for $0<\alpha\leq\frac 12$,
\be\label{1.1}
|\{(\ve_D: D); D\leq x, D \text{ nonsquare, } \ve_D \leq D^{\frac 12+\alpha}\}|\sim \frac {4\alpha^2}{\pi^2}x^{\frac 12}(\log x)^2.
\ee

In \cite {H}, Hooley also conjectures asymptotics for $\alpha >\frac 12$ of the form
\be\label{1.2}
|\{\cdots\}|\sim B(\alpha) x^{\frac 12} (\log x)^2
\ee
Thus in the range $\alpha_0<\alpha<c$, Theorem 1 recovers most of these discriminants.

\noindent
{\bf Remark 2.}
Exploiting Theorem 1.12 in \cite{M}, the Theorem will follow from the following statement.

{\sl $|\{D\leq x$; $D$ non-square and such that the Pell equation $t^2-Du^2=1$ has a solution with $t<x^{\frac 12+\alpha}$ 
and with the additional property
that for some integer $a<t, (a, t)=1$, $\frac at$ has partial quotients bounded by $A\}|$
\be\label{1.3}
\geq \big(1-o(1)\big) \frac{4\alpha^2}{\pi^{\frac 12}} x^{\frac 12} (\log x)^2.
\ee}

In order to prove \eqref{1.3}, we combine Hooley's approach recalled below with the results from \cite {B-K} around Zaremba's conjecture.
Denote
\be\label{1.4}
Z=\{t\in \mathbb Z_+; \text { there is $0<a<t, (a,t)=1$ such that $\frac at$ has partial quotients bounded by $A$}\}.
\ee

Observe that if $t^2-Du^2=1$ and $t<D$, then necessarily $t+u\sqrt D$ is the fundamental solution $\ve_D$, since otherwise
$$
2D> 2t> t+u\sqrt D\geq \ve_D^2\geq 4D
$$
using the property $\ve_D\geq 2\sqrt D$.
Hence, for $0<\alpha\leq \frac 12$, \eqref{1.3} will follow from the statement
\be\label{1.5}
\sum_{\substack{(t, u, D), t^2-Du^2=1\\ t\in Z, t< D^{\frac 12 +\alpha}, D< x}}
1> cx^{\frac 12} (\log x)^2 \text { with } c=\big(1-o(1)\big) \frac {4\alpha^2}{\pi^2}.
\ee
Denote
$$
\begin{aligned}
\mathcal R(u)= &\{\Omega(\mod u^2); \Omega^2\equiv 1 (\mod u^2)\}\\
&\rho(u) =|\mathcal R(u)|
\end{aligned}
$$
and let $X=\frac 12 x^\alpha$.
The left side of \eqref {1.5} is clearly minorized by
\be\label{1.6}
\sum_{1\leq u\leq X} \ \  \sum_{\Omega\in \mathcal R(u)} \  \sum_{\substack{ t\equiv\Omega (\mod u^2), t\in Z\\ u^{1+\frac 1{2\alpha}}< t< u\sqrt x}} 1
\quad\qquad\qquad\qquad\qquad\qquad\qquad\qquad
\ee
\begin{align}\label{1.7}
& =\sum_{1\leq u\leq X} \ \sum_{\Omega\in \mathcal R(u)} \ \  \sum_{\substack{ {u^{1+\frac 1{2\alpha}} <t< u\sqrt x} \\ {t\equiv\Omega (\mod
u^2)}}} 1 \\
&-\sum_{1\leq u\leq X} \rho(u)  \max_{\xi(\mod u^2)} |\{ t<u\sqrt x; t\not\in Z \text { and } t\equiv \xi
(\mod u^2)\}|.\label{1.8}
\end{align}

The expression \eqref{1.7} was evaluated in \cite{H} and \cite {F}.
We recall the argument. Since $u<\frac 12 x^\alpha$, clearly
\begin{align}\label{1.9}
\eqref{1.7}&\geq \sum_{1\leq u\leq X} \rho(u) \Big(\frac {u\sqrt x- u^{1+\frac 1{2\alpha}}} {u^2} -O(1)\Big)\nonumber\\
&>\big(1-o(1)\big) \Big\{\sum_{1\leq u\leq X}\frac {\rho (u)}u\sqrt x-\sum_{1\leq u\leq X} \rho(u) u^{\frac 1{2\alpha}-1}\Big\}.
\end{align}

Next, recalling \cite{F}, (24), (25)
\be\label{1.10}
\sum_{u\leq X} \frac {\rho(u)}u =\frac  4{\pi^2} (\log X)^2 +O(\log X)
\ee
\be\label{1.11}
\sum_{u\leq X}\rho (u)  u^{\frac 1{2\alpha}-1} <CX^{\frac 1{2\alpha}} \log X.
\ee
Substitution of \eqref{1.10}, \eqref{1.11} leads to the required minoration
$$
\big(1-o(1)\big) \frac 4{\pi^2} \sqrt x(\log X)^2 =\big(1-o(1)\big) \frac {4\alpha^2}{\pi^2}\sqrt x(\log x)^2.
$$
It remains to bound \eqref{1.8}.

We estimate \eqref{1.8} as follows
\begin{align}
\eqref {1.8} &\leq\sum_{\substack{{1\leq u\leq X}\\ {\rho(u)>(\log u)^{10}}}} \ \rho (u) \frac {\sqrt x}u\label{1.12}\\
&+\sum_{U\leq X} (\log U)^{10}\sum_{u\sim U} \max_{\xi(\mod u^2)} |\{ t\leq U\sqrt x; t\equiv \xi(\mod u^2), t\not\in Z\}|\label{1.13}
\end{align}
where $U$ takes dyadic values.

We rely on the following statement, to be proven in the next section.
\smallskip

\noindent
{\bf Proposition 2.}
{\sl There is a constant $c_1>0$ such that for $U<y^{c_1}$
\be\label{1.14}
\sum_{u\sim U} \max_\xi |\{t\leq y; t\equiv \xi (\text{\rm mod\,} u^2), t\not\in Z\}|< \frac yU(\log y)^{-20}
\ee
}

Estimate \eqref{1.12} by $O(\sqrt x)$ and, assuming
\be\label{1.15} 
\alpha<\frac {c_1}2
\ee
and using the Proposition, also
$$
\eqref{1.13} \lesssim \sum_{U\leq X} (\log U)^{10} \ \frac{U\sqrt x}U (\log x)^{-20} \lesssim \sqrt x.
$$
This proves the Theorem.

\smallskip
\noindent
{\bf Remark 3.} The proof of the Proposition appears in \S2 and depends essentially on he analysis in \cite{B-K} on Zaremba's
conjecture.
Let's point out here that the constant $c_1$ can be made explicit, since it relates to the `minor arcs analysis' in \cite{B-K}
and not to the spectral part.
Note that $|\{t<y; t\not\in Z\}|<y^{1-\delta}$ for some $\delta>0$ (see also \cite{M-O-W}), so that obviously (1.24) holds with $c_1<\frac \delta 2$.
The number $\delta$ depends however on the spectral theory of the continued fraction semigroup \break
$<\begin{pmatrix} a&1\\ 1&0\end{pmatrix} ; 1\leq a\leq A>^+$ and is non-explicit.

As mentioned above, our analysis permits to take $c=\frac 1{200}$ in Theorem 1, but as it stands certainly does not cover the full range $\alpha\in ]0, \frac 12]$ of Hooley's
result (1.1).

\section
{Proof of the Proposition}

We use the analysis from \cite {B-K} around the Zaremba problem.
In particular, the constant $A$ bounding the size of the partial quotients is taken at least $A\geq 51$ and may 
be further increased depending on our needs below.
Let $y=N$ be a large integer.
Following \cite {B-K}, replace $Z\cap [1, N]$ by a more convenient function $\theta: [1, N]\to [0, 1]$ with the property that $\supp
\theta\subset Z$.
This function is introduced in \cite {B-K} as the density on $Z\cap[1, N]$ obtained as a normalized
image measure of a suitable subset $\Omega_N$ of the
semi-group $\mathcal S_A$ under the map
$$
g\mapsto \langle g e_2, e_2 \rangle.
$$
A major part of the analysis in \cite {B-K} consisted then in the study of $\theta$ using the Hardy-Littlewood circle method and exploiting the
multilinear structure of $\Omega_N$.
Some of this analysis will also be used here.

We need to show that for $U<N^{c_1}$
\be\label{2.1}
\sum_{U\leq u < 2U}|\{n\leq N; n\equiv \xi_u (\mod u^2), n\not\in Z\}|<\frac NU (\log N)^{-20}
\ee
for any assignments $\xi_u(\mod u^2)$.

The condition $n\not\in Z$ will be replaced by $\theta(n)=0$.

The major arcs analysis in the circle method leads to a singular series density (cf. \cite{B-K}, \S4)
\be\label{2.2}
\sigma(n)\geq \prod_{p|n}\Big(1+\frac 1{p^2-1}\Big)\prod_{p|n} \Big(1-\frac 1{p+1}\Big)\geq \frac 1{\log\log n}.
\ee
Therefore it will suffice in order to establish \eqref{2.1} to show that
\be\label{2.3}
\sum_{U\leq u\leq 2U} \ \sum_{n\equiv\xi_u (\mod u^2)} |\theta(n) -\sigma(n)|^2 \lesssim \frac N{U(\log N)^{50}}.
\ee
Denote $\eta =\theta-\sigma$ and  $\hat\eta(\alpha) =\sum_n\eta(n) e(n\alpha), \alpha \in\mathbb R/\mathbb Z$, its Fourier transform.
By Parseval,
\begin{align}
\label{2.4}
\sum_{n\equiv\xi_u(\mod u^2)} |\eta(n)|^2 &= \frac 1{u^4} \int^1_0 \Big|\sum_{j(\mod u^2)} \hat\eta \Big(\frac j{u^2}+\alpha\Big)
e\Big(\frac {\xi _{u}}{u^2}j \Big)\Big|^2 d\alpha\nonumber\\
&\leq \frac 1{u^2} \int_{|\beta|<\frac 1{2u^2}} \Big[\sum_{j(\mod u^2)} \Big|\hat\eta \Big(\frac {j}{u^2}+\beta\Big)\Big|\Big]^2 d\beta
\end{align}
Define for dyadic $K$
$$
V_{q, K}=\{\alpha \in\mathbb T; \Big|\alpha-\frac aq\Big|\sim \frac KN \text{  for some } (a, q)=1\}
$$
setting
$$
V_{q, 1} =\{\alpha\in \mathbb T; \Big|\alpha-\frac aq\Big|\lesssim \frac 1N\text { for some } \ (a, q)=1\}.
$$
Since the function $\sigma$ was obtained in \cite{B-K} by restriction of $\hat\theta$ to the major arcs set
$$
\mathcal M=\bigcup_{\substack {q\leq Q_0\\ K< Q_0}} V_{q, K}
$$
with $Q_0 =e^{c(\log N)^{1/2}} $ say, we may in \eqref{2.4} replace $\hat\eta$ by $\hat\theta|_{\mathbb T\backslash \mathcal M}$.

Given $c_2< \frac 1{10^6}$ denote
\be\label{2.5}
V_{c_2}=\bigcup_{\substack {{q<N^{c_2}}\\ {K< N^{c_2}}}} V_{q, K}.
\ee
We recall several estimates from \cite{B-K} related to the minor arcs analysis.

\begin{lemma}\label{Lemma1}
\be\label{2.6}
\int_{\mathbb T\backslash V_{c_2}} |\hat\theta(\alpha)|^2 d\alpha < N^{1-c_3}
\ee
for some $c_3= c_3(c_2)>0$.
\end{lemma}

Denote
\be\label{2.7}
\hat\theta_0=\hat\theta 1_{V_{c_2}\backslash\mathcal M}.
\ee
From the preceding, we may estimate \eqref{2.4} by
\begin{align}
\label{2.8}
&\frac 1{U^2}\int_{|\beta|< \frac 1{2u^2}} \Big[\sum_{j(\mod u^2)}\Big|\hat\theta_0\Big(\frac j{u^2} +\beta\Big) \Big|\Big]^2
d\beta\\
&+\nonumber\\
&\sum_{j(\mod u^2)} \int_{|\beta|<\frac 1{U^2}} \Big|\hat\theta \Big(\frac j{u^2} +\beta\Big)\Big|^2 1_{\mathbb T\backslash V_{c_2}}
\Big(\frac j{u^2}+\beta\Big) d\beta.\label{2.9}
\end{align}

By Lemma 1, $\eqref{2.9}< N^{1-c_3}$ and the contribution in \eqref{2.3} at most $U.N^{1-c_3}$.
Choose in (2.5) the constant $c_2$ small enough to ensure the estimates below, borrowed from \cite {B-K}, valid.
Then (2,1) will hold for $c_1<\frac 12 c_3$.

In what follows, we always assume $q, K<N^{2c_2}$.

In the sequel, we use the notation $C$ for various numerical constants and let $A$ be large enough to make $C/A$ adequately small.

\begin{lemma}\label{Lemma2}
(see \cite{B-K}, Prop, 5.2 and Prop. 6.3).
\be\label{2.10}
|\hat\theta (\alpha)|\lesssim \frac N{(Kq)^{1-C/A}} \ \text { if } \alpha \in V_{q, K}.
\ee
\end{lemma}

The next statement is obtained by an easy variant of the proof of Prop. 6.21 in \cite{B-K}.

\begin{lemma}\label{Lemma3}
Assume given for each $q\sim Q$ a (possibly empty) subset $S_q$ of the residues $(\mod q)$.
Let $|\beta|< \frac KN$.
Then
\be\label{2.11}
\sum_{q\sim Q} \ \sum_{a\in S_q} \Big|\hat\theta \Big(\frac aq+\beta\Big)\Big| <\frac {(KQ)^{C/A}}{Q}
\Big(\sum_q |S_q|^2\Big)^{\frac 12} N.
\ee
\end{lemma}

A duality argument permits then to derive from \eqref{2.11}.

\begin{lemma}\label{Lemma4}
For $|\beta|<\frac KN$, we have
\be\label{2.12}
\Big(\sum_{\substack {{q\sim Q} \\ {(a, q)=1}}} \Big|\hat\theta\Big(\frac aq +\beta\Big)\Big|^2\Big)^{\frac 12} < \frac {(KQ)^{C/A} \log Q}{Q^{1/2}} N.
\ee
\end{lemma}
\smallskip

Hence 
\be\label{2.13}
\int_{\cup_{q\sim Q} V_{q, K}} |\hat\theta(\alpha)|^2 <\frac {(KQ)^{C/A}}{Q} KN.
\ee

Let $\Lambda \subset\{q\sim Q\}$.
Combining \eqref{2.13} with \eqref{2.10} gives

\begin{lemma}\label{Lemma5}
For $\Lambda$ as above, $Q<N^{2c_2}$
\be \label{2.14 }
\sum_{q\in \Lambda, (a, q)=1} \int_{|\beta|< \frac 1{N^{1-c_2}}} \Big|\hat\theta\Big(\frac aq
+\beta\Big)\Big|^2 d\beta < |\Lambda|^{\frac 12} \, \frac N{Q^{1-C/A}}.
\ee
\end{lemma}

We denote $(x, y)$ (resp. $[x, y]$) the $\ell cd$ (resp. $s cm$) of $x, y\in \mathbb Z^*$.
Returning to (2.8) fix $u\sim U$ and assume $\frac {j_0}{u^2} +\beta \in V_{q, K}$ with $q, K<N^{c_2}$  for some $j_0(\mod u^2)$. Hence
\be\label{2.15}
\frac {j_0}{u^2} +\beta =\frac {a_0}{q_0}+\gamma 
\ee
where $q_0< N^{c_2}, (a_0, q_0)=1$ if $a_0\not= 0$ and $|\gamma|<\frac 1{N^{1-c_2}}$.
Then for all $j(\mod u^2)$
$$
\frac j{u^2} +\beta =\frac {j-j_0}{u^2}+\frac {a_0}{q_0}+\gamma \in V_{q', N^{c_2}}\ \text {  where  } \ q'|u^2 q_0, u^2q_0< N^{2c_2}.
$$

In particular
\be\label {2.16}
\beta = - \frac {j_0}{u^2}+\frac {a_0}{q_0} +\gamma =\beta_1+\gamma.
\ee

Denote $\nu_p {(x)}$ the exponent of the prime $p$ in the factorization of $x\in\mathbb Z^*$.
Decompose $\beta_1$ as
\be\label{2.17}
\beta_1 =\frac {a_1}{q_1} +\frac {a_2}{q_2}+\frac {a_3}{q_3} =\frac {a_1}{q_1}+\tilde\beta_1
\ee
where $(a_i, q_i)=1$ if $a_i\not= 0$ and
\begin{itemize}
\item[(2.18)] \  $q_1|u^2$
\item
[(2.19)] \ If $\nu_p (q_2)>0$, then $\nu_p(q_2)> 2\nu_p(u)>0$ and $\big|\frac {a_2}{q_2}\big| < \frac 1{(q_2, u^2)}$
\item [(2.20)] \ $(a, q_3)= 1$.
\end{itemize}

Note that since $|\beta_1|\leq |\beta| +|\gamma|< \frac 1{U^2} +\frac 1{N^{1-c_2}} <\frac 2{U^2}$ and $\beta_1-\tilde\beta_1\in \frac 1{u^2} \mathbb Z$ by (2.17),
$\tilde \beta_1$ essentially determines $\beta_1$. Next
$\frac j{u^2}+\beta_1\in \big\{ \frac {j'}{u^2}; j' (\mod u^2)\big\} +\tilde\beta_1$, hence $\frac j{u^2}+\beta_1=\frac aq +\tilde\beta_1 =\frac {a'}{q'}$ where
$q|u^2, (a, q)=1= (a', q')$.
Note that $\tilde \beta_1$ belongs to the set $\mathfrak S_q$ of sums $\frac {a_2}{q_2} +\frac {a_3}{q_3}$ where $(a_i, q_i)= 1$ if $a_i \not=0$ and
\begin{itemize}
\item [(2.21)] \ If $\nu_p(q_2)>0$, then $\nu_p(q_2)> \nu_p(q)$ and $\big|\frac {a_2}{q_2}\big|<\frac 1{(q_2, q)}$
\item [(2.22)] \  $(q_3, qq_2)=1$
\end{itemize}
as a consequence of (2.19), (2.20).
With previous notation, we have for $\tilde\beta_1 \in \mathfrak S_q$ that
$$
\frac {a'}{q'} =\frac aq +\frac {a_2}{q_2} +\frac {a_3}{q_3} =  \frac {\tilde a}{[q, q_2]}+\frac {a_3}{q_3} \ \text { with } \  \tilde a =
\frac {q_2}{(q, q_2)} a+ \frac {q}{(q, q_2)} a_2, (\tilde a, qq_2)=1
$$
and $q' =[q, q_2]q_3$.
We claim that $\frac {a'}{q'}$ essentially determines $\frac aq, \tilde\beta_1$.
Fix $q_1, q_2, q_3$ divisors of $q'$.
Since $a' \equiv a_3 [q, q_2](\mod q_3), ([q, q_2] q_3)= 1, a_3$ is determined and hence $\tilde a$.
Recalling (2.21), $|a_2|< \frac {q_2}{(q, q_2)}$ and since $\frac q{(q, q_2)} a_2\equiv \tilde a\big(\mod  \frac {q_2} {(q, q_2)}\big)$, $a_2$ is determined.
This proves the claim.

Estimate using Cauchy-Schwarz inequality
$$
\sum_{u\sim U} (2.8) \leq \frac 1{U^2} \sum_q q\Big(\sum_{q|u^2, u\sim U} d(u^2)\Big) \sum_{\tilde\beta_1 \in\mathfrak S_q}
\sum_{(a, q)=1}\int_{|\gamma|< \frac 1{N^{1-c_2}}} \Big|\hat\theta_0 \Big(\frac aq +\tilde\beta_1+\gamma\Big)\Big|^2
d\gamma.\eqno{(2.23)}
$$

If $q|u^2$, there is $q_1|u$ such that $q_1|q, q|q_1^2$. We obtain
$$
\begin{aligned}
(2.23) &\leq \frac 1{U^2}\sum_{q_1\leq U} d(q_1)^2 \sum_{q_1|u, u\sim U} d\Big(\frac u{q_1}\Big)^2 \sum_{q_1|q, q|q_1^2} q \, 
\sum_{\substack{\tilde \beta_1\in\mathfrak S_q\\ (a, q)=1}} \int_{|\gamma|<\frac 1{N^{1-c_2}}} |\hat\theta_0 \Big(\frac aq +\tilde\beta_1+\gamma\Big)|^2 d\gamma\\
&\lesssim \frac {(\log U)^{10}}U \ \sum_{q_1\leq U} \ \frac {d(q_1)^2}{q_1}\sum_{\substack{q_1|q\\ q|q_1^2}} q \sum_{q|q', q'<N^{2c_2}}
\sum_{(a', q')=1} \ \int _{|\gamma|<\frac 1{N^{1-c_2}}} \ \Big|\hat\theta_0\Big(\frac {a'}{q'}+\gamma\Big)\Big|^2 d\gamma
\end{aligned}
\eqno{(2.24)}
$$
in view of previous discussion.

Define for $Q_1\leq Q\leq Q^2_1$, $Q\leq Q'< N^{2c_2}$ ranging dyadically, the set 
$$
\Lambda_{Q', Q, Q_1} =\big\{q' \sim Q'; \text { there are $q_1\sim Q_1, q\sim Q$ such that $q_1|q, q|q_1^2, q|q'\big\}$}.
$$
Hence
$$
|\Lambda_{Q', Q, Q_1}|\ll Q_1^{1+} \frac {Q'}Q.\eqno{(2.25)}
$$
Restricting $q_1\sim Q_1, q\sim Q, q'\sim Q'$, using Lemma 5, the corresponding contribution to (2.24) may be bounded by
$$
\begin{aligned}
\frac {(\log U)^{10}}U \ \frac Q{Q_1^{1-}} \big|\Lambda_{Q', Q, Q_1}\big|^{\frac 12} \ \frac N{(Q')^{1-C/A}}&\ll
\frac {(\log U)^{10}}{U} \, Q_1^{-\frac 12} Q^{\frac 12} (Q')^{-\frac 12+C/A} N\\
&\ll \frac {(\log U)^{10}}{U} \, Q_1^{-\frac 12+\frac {4C}{A}} (Q')^{-C/A} N.
\end{aligned} \eqno{(2.26)}
$$
Assuming $A$ large enough and recalling  (2.7) restricting $Q'>Q_0$, leads to the estimate
$$
\frac {(\log U)^{10}}{U} \ Q_0^{-C/A} <\frac NU(\log N)^{-20}.\eqno{(2.27)}
$$
This proves (2.1) and Proposition 2.

\end{document}